\documentclass[11pt]{article}

\usepackage[T1]{fontenc}
\usepackage[latin1]{inputenc}
\usepackage[cm]{aeguill}

\usepackage{amsmath}
\usepackage{amssymb}
\usepackage{amsthm}
\usepackage{nicefrac}
\usepackage{graphicx}

\theoremstyle{plain}
\newtheorem{theorem}{Theorem}[section]
\newtheorem{proposition}[theorem]{Proposition}
\newtheorem{lemma}[theorem]{Lemma}
\newtheorem{remark}[theorem]{Remark}
\newtheorem{assumption}[theorem]{Assumption}
\newtheorem{corollary}[theorem]{Corollary}

\theoremstyle{definition}
\newtheorem{proc}[theorem]{Tuning procedure}

\newcommand{\eps}{\varepsilon}                        
\renewcommand{\phi}{\varphi}                            

\newcommand{\nr}[1]{\left\Vert #1\right\Vert}         
\newcommand{\abs}[1]{\left\vert #1\right\vert}        

\DeclareMathAlphabet{\mathonebb}{U}{bbold}{m}{n} 
\newcommand{\one}{\ensuremath{\mathbf{1}}}

\newcommand{\sur}[3][.5ex]{%
  \raise#1\hbox{\ensuremath{#2}}\!\!\left/%
    \vphantom{\raise#1\hbox{\ensuremath{#2}}}\!%
    \lower#1\hbox{\ensuremath{#3}} \right.%
}

                    % ...

\newcommand{\E}{\mathbb{E}}                          
\newcommand{\R}{\mathbb{R}}                          
\newcommand{\N}{\mathbb{N}}                          
\newcommand{\Pb}{\mathbb{P}}                          

\title{Simulation of a Local Time Fractional Stable Motion}
\author{Matthieu Marouby \medskip
\\
\small\textit{Institut de Mathématiques, Université Paul Sabatier,} \\
\small\textit{31062 Toulouse, France}\\
\small\textit{email address : marouby@math.univ-toulouse.fr}}
\date{}

\begin{document}
\bibliographystyle{plain}

\maketitle

\abstract
In this paper, we simulate sample paths of a class of symmetric $\alpha$-stable processes using their series
expression. We will develop a result in the approximation of shot-noise series. And finally, we will get a convergence rate for the approximation.
\\

\\


\section{Introduction}
Fractional fields have often been used to model irregular phenomena. The simplest one is the fractional Brownian motion introduced in \cite{kolmogorov1940wse} then developed in \cite{mandelbrot1968fbm}. 
More recently, many fractional processes have been studied, usually obtained by a stochastic integration of a
deterministic kernel against a random measure (cf. among others \cite{levy2000rrp},
\cite{benassi2002iap}, \cite{kaj2007cfb}, \cite{lacaux2004rhm} and \cite{MR2290879}). 
Many different simulations methods were discussed in the literature, but shot noise series seem to fit perfectly in that kind of problem. Generalized shot noise series
where introduced for simulation in \cite{rosinski1987ppc}, further developments were done in \cite{rosinski} and \cite{MR1833707} and a general framework was developed in
\cite{cohen2007gfs}. Moreover, a computer study of the convergence rate of LePage series to $\alpha$-stable random variables has been done in a particular case in \cite{janicki1992cir}.

A shot noise series can be seen as: $$\sum_{n=1}^{\infty} \Gamma_{n}^{-1/\alpha} V_{n}$$
where $(V_{n})$ are i.i.d. random variables and $(\Gamma_{n})$ the arrival times of a Poisson process. Usually, there
is no question about the simulation of $V_{n}$. In this paper, we will contribute to this theory by looking at the
convergence rate when $V_{n}$ can not be simulated but only approximated. Moreover, we will show a way 
to save computer time. Indeed, in a shot noise series representation of a $\alpha$-stable process,
the first terms are the bigger so you have to minimize the error made in approximating $V_{n}$ for small $n$. For large
$n$, $ \Gamma_{n}^{-1/\alpha}$ is small that it is not as useful to approximate $V_{n}$
them with so much care. We developed a way of measuring the convergence rate towards the limiting process depending
on the error made while computing each term.

Afterwards, we will consider an application of this result to study a particular process. In network
traffic modeling, properties like self-similarity, heavy tails and long-range dependance are often needed, see for
example \cite{mikosch2002nta}. 
Moreover, empirical studies like \cite{crovella1997ssw} showed the importance of self-similarity and long-range
dependance in that area.

In \cite{cohensamo}, the authors introduced fractional Brownian motion local time fractional stable motion as a stochastic integration of a non-deterministic
kernel against a random measure, which will be our main interest in the second part. Here we will call it Local
Time Fractional Stable Motion (LTFSM).
This process was defined as:
$$\int_{\Omega'}\int_{\R} l(x,t)(\omega') M(d\omega',dx), \text{ for } t\geq 0.$$
In this expression, 
$l$ is the local time of a fractional Brownian motion of Hurst parameter $H$ defined on $(\Omega',\mathcal{F}',\mathbf{P'})$. $M$ is a symmetric alpha stable
random measure (see \cite{taqqusamo} for more details) with control measure $\mathbf{P}'\times Leb$ (Leb being the
Lebesgue measure on $\R$). 
LTFSM is $\alpha$-stable but also self-similar and its increments are long-range dependant.

The first step towards understanding LTFSM is naturally to have a look at its sample paths. Unfortunately, the
above expression does not directly give a way to get the sample paths. In the case of Brownian motion local time,
this process can be seen as the limit of a discrete random walk with random rewards model. It is not completely
satisfying for a few reasons: first, it only works for $H=1/2$, then, there is no control of the convergence speed rate
towards the limit. That is where the tool we developed in the first part comes to play. In this paper, we will study how
we can simulate this process by using the expression given in equation (5.3) in \cite{cohensamo}, which can be seen as
a shot noise expansion.

The next section is devoted to the shot noise theory results.
In the following part we will see how LTFSM fits in the general frame we just developed,
and we will be able to get a convergence rate in the case of our process.
The last section will be devoted to a quick study of our simulations with a comparison to the random walk with random rewards model that was introduced in \cite{cohensamo} in the case
$H=1/2$.

\section{Shot Noise series}\label{shotnoise}
In this section, we will show some results on shot noise series, using theorem 2.4 in \cite{rosinski}. 
\begin{assumption}\label{ass1}
$E_{K}$ will be the space of continuous functions defined on a compact subset $K \subset \R$, equipped with the uniform norm denoted $\nr{\cdot}_{K}$. We will denote by $\nr{\cdot}_{K,p}$ the $L^p(K)$ norm.
Let us consider the case $h : \R \times E_{K} \rightarrow E_{K}$ where  
\begin{equation}\label{defh}
h(r,v)= r^{-1/\alpha} v.
\end{equation}
$h$ is a Borel measurable map. 
Let $(\Gamma_{n})_{n\geq 1}$ be the arrival times of a Poisson process of rate $1$ in $\R_{+}$, $(V_{n})_{n \geq 1}$ be a sequence of i.i.d. symmetrical random variables of distribution $\lambda$. Let us assume $(\Gamma_{n})_{n \in \N}$, $(V_{n})_{n \in \N}$ are independent.
\end{assumption}
We will suppose that the assumptions of corollary 2.5 in \cite{rosinski} are satisfied, i.e. :
\begin{assumption}\label{ass2}
The measure defined for all $A \in \mathcal{B}_{E}$ (Borelians of $E_{K}$) by 
\begin{equation}\label{deflevymeasure}
F(A) = \int_{0}^{\infty} \int_{E_{K}} \one_{A\setminus \{0\}} (h(u,v)) \lambda (dv)du
\end{equation}
is a Lévy measure.
\end{assumption}

Now, we will get some inspiration from the proof of \cite{lacaux} to prove the following proposition, considering the $L^q$ distance between the sum of the series, and the truncated sum.

\begin{proposition}\label{evalreste}
Under assumptions \ref{ass1} and \ref{ass2}, let us consider the shot noise series $$Y=\sum_{n=1}^{\infty} h(\Gamma_{n}, V_{n}).$$
Let us denote $Y_{N}= \sum_{n=1}^{N} h(\Gamma_{n}, V_{n}).$
Take $q \geq 2$, we have to suppose that there exists $M_{q}$ such that for all $n$, $\E[\nr{V_{n}}_{K}^{q}] \leq M_{q} < \infty$.
Then, for $N > q / \alpha - 1$, we have
$$\E[\abs{Y(t)-Y_{N}(t)}^{q}] \leq \frac{A_{q} H_{N+1,q}}{N^{q(\frac{2- \alpha}{2 \alpha})}},$$
where $A_{q}$ depends only on $q$ through $M_{q}$ and $B_{q}$, and 
$$H_{n,q} = \frac{\Gamma\left(n  - \frac{q}{\alpha} \right) n^{q/\alpha}}{\Gamma(n)}.$$
with $\lim_{n \rightarrow \infty} H_{n,q}=1$.
\end{proposition}

\begin{proof}
Denote for $N<P$ : $$R_{N,P}(t) = \sum_{n=N+1}^{P} h(\Gamma_{n},V_{n})(t).$$
According to theorem 2.4 in \cite{rosinski}, $\lim_{P \rightarrow + \infty} R_{N,P}$ exists.

As $h(\Gamma_{n},V_{n})(t)= \Gamma_{n}^{-1/\alpha} V_{n}(t)$ and because $V_{n}(t)$ is symmetric, $h(\Gamma_{n},V_{n})(t)$
is also symmetric. Therefore, we can apply proposition 2.3 in \cite{ledoux1991pbs}, and we get
$$\E \left[ \max_{N+1 \leq n \leq P} \abs{R_{N,n}(t)}^{q} \right] \leq 2 \E[\abs{R_{N,P}(t)}^{q}].$$ 
We now are able to apply Khintchine inequality
Let $\eps_{n}$ be a sequence of i.i.d. Rademacher random variables, independent from everything else. Thus $\eps_{n} h(\Gamma_{n},V_{n})(t)$ has the same distribution as $h(\Gamma_{n},V_{n})(t)$ since it is symmetrical.
Khintchine inequality claims
 \begin{multline}
 \E\left[\abs{\sum_{n=N+1}^{P} \eps_{n} h(\Gamma_{n},V_{n})(t) }^{q} \Biggl\vert (h(\Gamma_{n},V_{n})(t))_{n \in \N} \right]^{1/q} \leq \\
  B_{q} \left( \sum_{n=N+1}^{P} \abs{h(\Gamma_{n},V_{n})(t)}^{2} \right)^{1/2},$$
\nonumber\end{multline}
where $B_{q} = \sqrt{2} \left(\frac{\Gamma((q+1)/2)}{\sqrt{\pi}} \right)^{1/q}$ for $q>2$, and $B_{q}= 1$ if $q=2$.
%
Taking the expected value on both sides of the inequality, then using Minkowski's inequality, we have
\begin{equation}\label{major1}
\E\left[\abs{\sum_{n=N+1}^{P}h(\Gamma_{n},V_{n})(t)}^{q} \right] \leq B_{q}^{q} \left( \sum_{n=N+1}^{P} \E\left[\abs{h(\Gamma_{n},V_{n})(t)}^{q}\right]^{2/q}\right)^{q/2}.
\end{equation}

We now, have only to compute  $\E\left[ \abs{h(\Gamma_{n},V_{n})(t)}^{q}\right]$.
As $V_{n}$ and $\Gamma_{n}$ are independent, as we can compute $\E[\Gamma_{n}^{-q/\alpha}]$ and $\E[\nr{V_{n}}_{K}^{q}]
\leq M_{q} <\infty$, we have
$$\E\left[ \abs{h(\Gamma_{n},V_{n})(t)}^{q}\right] \leq M_{q} \frac{\Gamma\left(n- \frac{q}{\alpha} \right)}{\Gamma(n)} .$$
We can use this results in (\ref{major1}), leading to
\begin{equation}\label{majorreste}
\E[\abs{R_{N,P}(t)}^{q}] \leq  B_{q}^{q} M_{q} \left( \sum_{n=N+1}^{P}\left( \frac{\Gamma\left(n- \frac{q}{\alpha} \right)}{\Gamma(n)} \right)^{2/q} \right)^{q/2}.
\end{equation}
Stirling formulae claims
$$\sum_{n=N+1}^{+\infty}\left( \frac{\Gamma\left(n- \frac{q}{\alpha} \right)}{\Gamma(n)} \right)^{2/q} \sim \sum_{n=N+1}^{+\infty} \frac{1}{n^{2/\alpha}} \sim \frac{\alpha}{(2-\alpha) N^{2/\alpha - 1}}.$$
which shows us that the series converges. Let us denote
\begin{equation}\label{defhnq}
H_{n,q} = \frac{\Gamma\left(n  - \frac{q}{\alpha} \right) n^{q/\alpha}}{\Gamma(n)}.
\end{equation}
It can be easily proved that $\sup_{n \geq N+1} H_{n,q} = H_{N+1,q}$, and we have
$$\E \left[ \sup_{M \geq N+1} \abs{R_{N,M}(t)}^{q} \right] \leq \frac{A_{q} H_{N+1,q}}{N^{q/2 (2/ \alpha - 1)}},$$ 
with $$A_{q}= 2 B_{q}^{q} M_{q} \left( \frac{\alpha}{2-\alpha} \right)^{q/2}.$$
\end{proof}

We can improve this result by checking what happens on $E_{K}$ equipped with $\nr{\cdot}_{K,p}$. 
\begin{proposition}\label{erreurlplq}
For $p>0$ and $(N+1) \alpha >q > \max(p,2)$, there exists $A_{q}$ such that
$$\E \left[ \nr{Y- Y_{N}}_{K,p}^{q} \right] \leq Vol(K)^{q/p} A_{q} \frac{H_{N,q}}{N^{q (\frac{2-\alpha}{2\alpha})}}$$
\end{proposition}

\begin{proof}
According to Hölder's inequality, we have $$\nr{Y-Y_{N}}_{K,p} \leq \left( Vol(K)^{1-p/q}  \left( \int_{K} \abs{Y(t)-Y_{N}(t)}^{q}dt \right)^{p/q}   \right)^{1/p},$$
where $Vol(K)=\int_{\R} \one_{K}$.
Thanks to  proposition \ref{evalreste} we have
$$\E \left[ \nr{Y- Y_{N}}_{K,p}^{q} \right] \leq Vol(K)^{q/p-1} \int_{K} \frac{A_{q} H_{N,q}}{N^{q/2 (2/\alpha - 1)}}dt.$$
\end{proof}

We have proved that we are close enough to the process $Y$ if we truncate the sum. Unfortunately, we are not always able to simulate $V_{n}$. We will now try to see what happens if we use a sequence of random variables $(W_{n,k})_{n \geq 1}$ such that $\lim_{k \rightarrow \infty} W_{n,k} = V_{n}$ in a sense we will define later.
In the next proposition, we will evaluate the distance between $\sum_{n=N+1}^{P} h(\Gamma_{n},V_{n})(t)$ and $\sum_{n=N+1}^{P} h(\Gamma_{n},W_{n,k})(t)$ in $L^q$ for $q >2$.
Due to the fact that $\Gamma_{n}^{-1/ \alpha}$ has a q-moment if and only if $n>q / \alpha$, we will not be able to compute in general the distance between the sums starting at $n=1$.

\begin{proposition}\label{evalsomme}
For $\alpha(P+1)>\alpha(N+1)>q\geq 2$, if $(W_{n,k})_{n\geq 1}$ is a sequence of random variables such that there
exists a constant $M_{q,k}$: 
$$\E[\nr{V_{n}-W_{n,k}}_{K}^{q}] \leq M_{q,k} n^{q \beta} < \infty,$$
with $\beta < 1/\alpha -1/2$.
Then
\begin{multline}
$$\E \left[ \abs{\sum_{n=N+1}^{P}h(\Gamma_{n},V_{n})(t) - h(\Gamma_{n}, W_{n,k})(t)}^{q} \right] \leq   A_{q}' H_{N+1,q} M_{q,k} \\
\left( \frac{1}{N^{2/\alpha-\beta-1}} - \frac{1}{P^{2/\alpha-\beta-1}} \right)^{\frac{q}{2}} \nonumber
\end{multline}
where $A_{q}'$ depends only on $q$.
\end{proposition}
%

\begin{proof}
Let us denote  

$$R_{N,P}(t)= \sum_{n=N+1}^{P} h(\Gamma_{n},V_{n})(t) - h(\Gamma_{n},W_{n,k})(t).$$

In the same fashion we got (\ref{majorreste}) in the proof of proposition \ref{evalreste}, we have
$$\E[\abs{R_{N,P}(t)}^{q}] \leq  B_{q}^{q} M_{q,k} \left( \sum_{n=N+1}^{P}\left( n^{q \beta} \frac{\Gamma\left(n- \frac{q}{\alpha} \right)}{\Gamma(n)} \right)^{2/q} \right)^{q/2}.$$

Using the same definition of $H_{n,q}$ as we did in (\ref{defhnq}), and knowing that $\sup_{n \geq N+1} H_{n,q} = H_{N+1,q}$,
and using well known series-integral comparison we get
$$\E[\abs{R_{N,P}(t)}^{q}] \leq   A_{q}' H_{N+1,q} M_{q,k} \left( \frac{1}{N^{2/\alpha - \beta-1}} - \frac{1}{P^{2/\alpha - \beta -1}} \right)^{q/2}$$
where $A_{q}'= B_{q}^{q} \left( \frac{\alpha}{2 - \alpha \beta-\alpha} \right)^{q/2}$.

\end{proof}
In the same way we had proposition $\ref{erreurlplq}$, we have
\begin{proposition}\label{evalsommelplq}
For $p>0$ and $(P+1)\alpha  > (N+1) \alpha >q > \max(p,2)$, under the same assumptions as in proposition \ref{evalsomme} we have
\begin{multline}
$$\E \left[ \nr{\sum_{n=N+1}^{P}h(\Gamma_{n},V_{n}) - h(\Gamma_{n}, W_{n,k})}_{K,p}^{q} \right] \leq Vol(K)^{q/p} A_{q}'
\\ H_{N+1,q} M_{q,k} \left( \frac{1}{N^{2/\alpha-\beta-1}} - \frac{1}{P^{2/\alpha-\beta-1}} \right)^{\frac{q}{2}}.
\nonumber
\end{multline}
\end{proposition}

The previous proposition let us evaluate the error we do approximating in the sum $V_{n}(t)$ by $W_{n,k}(t)$, knowing their $L^q$-distance. In the next part, we will use these results by taking enough terms to minimise the error done in proposition \ref{evalreste} and a good enough approximation of $V_{n}$ such that $M_{q,k}$ is small, in order to have a small error computing $\sum h(\Gamma_{n},W_{n,k})$ instead of $\sum h(\Gamma_{n},V_{n})$.
Depending on the values of $q$ and $\alpha$, we may have to deal with the first terms in a particular way.


\section{Application to the local time fractional stable motion}
In this section, we will apply the results from section \ref{shotnoise} on the process defined in the introduction. The precise definition we will work on is:

\begin{equation}\label{defy}
\sum_{n \in \N}  \Gamma_{n}^{-1 / \alpha} G_{n}  e^{X'^{2}_{n}/ 2 \alpha} l_{n}(X'_{n},t),
\end{equation}
where
\begin{itemize}
    \item $(G_{n})_{n \geq 1}$ is a sequence of i.i.d. of standard normal distribution,  
    \item $(X'_{n})_{n \geq 1}$ is another one,
    \item $(l_{n})_{n \geq 1}$ are independent copies of a fractional Brownian motion local time, each one defined on some $(\Omega',\mathcal{F}',P')$. 
    \item $(\Gamma_{n})_{n \geq 1}$ are the arrival time of a Poisson process of rate $1$ on $[0,\infty)$,
\end{itemize}
all the variables being independent.

But before that, we will have to prove that this process satisfies the assumptions needed.

In fact, we will be able to generalise a bit this work. Consider $f_{n}$ satisfying the following assumptions: 
$f_{n}: \R \times \R_{+} \rightarrow \R$ are the independent copies of a non negative continuous random function on the probability space $(\Omega_{n}',\mathcal{F}_{n}',P_{n}')$,
such that for all $K \subset \R$ compact set, denoting $E_{K}$ the space of continuous functions on $K$ equipped with the uniform norm $\nr{\cdot}_{K}$, we have for $p>0$

$$\E \left[  \nr{f_{n}(x,\cdot)}_{K}^{p} \right] $$ uniformly bounded in $x$, and $f_{n}(\cdot,t)$ has its support included in
\begin{equation}\label{defsrho}
S_{\rho,n}= \left[\inf_{s \leq t} B_{s}^{H,n} - \rho,  \sup_{s \leq t} B_{s}^{H,n} + \rho \right],
\end{equation}
 where $(B^{H,n})$ are independent fractional Brownian motion with Hurst parameter $H$ defined on the same space as $f_{n}$ and $\rho \geq 0$. 
In the following, we will consider 
\begin{equation}\label{defy}
Y(t) = \sum_{n \in \N} \Gamma_{n}^{-1 / \alpha} G_{n}  e^{X'^{2}_{n}/ 2 \alpha} f_{n}(X'_{n},t),
\end{equation}

\begin{remark}
The local time of a fractional Brownian motion satisfies obviously the support condition with $\rho=0$.
It satisfies the other condition because $l_{n}(x,\cdot)$ is a non decreasing function so that $$\sup_{t\in K} l_{n}(x,t)=l_{n}(x,t_{0}),$$
 $t_{0}$ being the upper bound of $K$. It simply claims $\E \left[ \nr{l_{n}(x,t)}_{K}^{p}\right] \leq t_{0}^{p}$.
\end{remark}


The following lemma is a direct consequence of the support condition so we will skip its proof.

\begin{lemma}\label{lemmeholder}
Let $f_{n}(x,\cdot)$ be a continuous function on K, $\nr{\cdot}_{K}$ being the uniform norm on $K$.
Let $f_{n}$ satisfies the assumptions stated above. Let $\phi$ be the positive density of $X_{n}$ distribution with respect to Lebesgue measure. For $q>0$, if for all $a>0$, $\int_{\R} \phi(x)^{-q/\alpha+1} e^{- a x^{2}} < \infty$,
there exists $C$ such that 
$$\E \left[\nr{\phi(X_{n})^{-1/ \alpha} f_{n}(X_{n},\cdot)}_{K}^{q} \right] \leq C \sup_{x\in \R}\left( \E\left[\nr{f_{n}(x,\cdot)}_{K}^{qp'}\right]\right)^{1/p'},$$
for some $p'>1$.
\end{lemma}

%

For some technical reasons that will appear later, we would like to work with a slightly different process, having the same distribution as $Y$. Because of that, we will show that the process $Y$ and $Y_{\phi}$ are identically distributed, where
$$Y_{\phi}(t) = \sum_{n \in \N} \Gamma_{n}^{-1 / \alpha}  G_{n} \phi(X_{n})^{-1/ \alpha} f_{n}(X_{n},t),$$
where $X_{n}$ are i.i.d. random variables, and $\phi$ is the density of their distribution which satisfies for $q>0$, for all $a>0$, $\int_{\R} \phi(x)^{-q/\alpha+1} e^{- a x^{2}} < \infty$.

In order to achieve that objective, we will show that $Y$ and $Y_{\phi}$ have the same characteristic function. 
We can compute this characteristic function by using techniques that can be found in \cite{lacaux}, thus we get
\begin{proposition}\label{carac}
Let $K$ be a compact subset of  $\R_{+}$, $E_{K}$ be the space of continuous functions on $K$,  equipped with the uniform norm $\nr{\cdot}$, and $E'_{K}$ its dual space. For all $y' \in E_{K}'$, we have
\begin{multline}
\E \left[ e^{i <y', Y_{\phi}(\cdot) >} \right] =\\
\exp \int \left( e^{i < y', h(r,v)>} - 1 - i <y', h(r,v)> \one_{\nr{h(r,v)} \leq 1}  \right) dr \lambda(dv), \nonumber$$
\end{multline}
where $\lambda$ is the distribution of $V_{n} = G_{n} \phi(X_{n})^{-1/ \alpha} f_{n}(X_{n}, \cdot)$ and
\begin{align}
h : \R_{+} \times E_{K}& \rightarrow E_{K} \nonumber\\
(t,v(\cdot))& \mapsto t^{-1/\alpha} v(\cdot). \nonumber
\end{align}
\end{proposition}

Then straightforward calculus lead to the following proposition.

\begin{proposition}\label{samedistrib}
If for some $q>0$ and for all $a>0$ we have $\int_{\R} \phi(x)^{-q/\alpha+1} e^{- a x^{2}} < \infty$,
then processes $Y_{\phi}$ and $Y$ have the same distribution.
\end{proposition}

From now on, we will only use 
\begin{equation}\label{defy2}
Y(\cdot) = \sum_{n \geq 1} \Gamma_{n}^{-1 / \alpha} G_{n}  e^{ 2 \abs{X_{n}}/ \alpha} f_{n}(X_{n},\cdot),
\end{equation}
where $X_{n}$ has a Laplace distribution of parameters $(0,1/2)$, i.e. its density is $e^{-2 \abs{x}}$ with respect to the Lebesgue measure on $\R$.

\begin{proposition}\label{cvuni}
Let $K \subset \R_{+}$ denote a compact set. The series defining the process $Y(t)$ in (\ref{defy2}) converges uniformly on $K$.
\end{proposition}

\begin{proof}
The proof of this proposition simply consists in verifying that $Y$ satisfies the assumptions of theorem 2.4 in \cite{rosinski}, that
is:
The measure $F$ defined for all borelians in $E_{K}$ in (\ref{deflevymeasure}) is a Lévy measure. We have to prove
\begin{itemize}
  \item[(i)]  $$\int_{E_{K}} (<y',y>^{2} \wedge 1)  F(dy) < \infty$$
  \item[(ii)] The function 
        \begin{equation}
        \Phi: y' \mapsto \exp \left( \int_{E_{K}} \left( e^{i <y',y>} -1 -i <y',y> \one_{\nr{y} \leq 1} \right) F(dy)
\right)\nonumber
        \end{equation}
        is the characteristic function of a probability measure on $E_{K}$.
\end{itemize}

(ii) is a direct consequence of proposition \ref{carac}.
To prove (i), a straightforward computation leads to: 
$$\int_{E_{K}} (\nr{y}^{2} \wedge 1)  F(dy) = \int_{\R} \abs{v}^{\alpha} \lambda'(dv) \left( 1 +  \frac{\alpha}{2 -
\alpha} \right).$$
where the second member of the equality is finite because by independence of the $G_{n}$: 
$$\int_{\R} \abs{v}^{\alpha} \lambda'(dv) = \E[\abs{G_{n}}^{\alpha}] \E[\phi(X_{n})^{-1} \sup_{t \in K} f_{n}(X_{n},t)^{\alpha}].$$
Using lemma \ref{lemmeholder} the integral is finite.
Theorem 2.4 in \cite{rosinski} can thus be applied.

\end{proof}

We will now see how we can apply propositions \ref{erreurlplq} and \ref{evalsommelplq}


\begin{remark}
For $K \subset \R_{+}$ compact, 
$$M_{q} := \E [\nr{V_{n}}_{K}^{q}] < \infty. $$

Indeed, by independence of the variables,
$$M_{q} = \E[\abs{G_{n}}^{q}]  \E\left[e^{\frac{2 \abs{X_{n}}q}{\alpha}} \nr{f_{n}(X_{n},t)}_{K}^{q}\right],$$
and the last expectation can be bounded using lemma \ref{lemmeholder}.
\end{remark}

Using proposition \ref{erreurlplq},  we directly get the following corollary
\begin{corollary}\label{erreurreste}
For $K \subset \R_{+}$ compact, there is a constant $C$ such that for $p>0$ and $P (\alpha +1)>q >\max(p,2)$
$$\E \left[\nr{Y-Y_{P}}_{K,p}^{q}\right] \leq \frac{C}{P^{q(\frac{2- \alpha}{2 \alpha})}}$$
\end{corollary}

Now, let us study the non-truncated terms. Denoting $$g_{n,k}(x,t)= \int_{\R} \phi_{k}(y-x)l_{n}(y,t)dy,$$ where $(\phi_{k})_{k}$ is an approximate identity with support in $[-1/k,1/k]$. 
We will use $\phi(x)= -\abs{x}+1$ on $[-1,1]$, $\phi=0$ elsewhere. We will denote $\phi_{k}(x)= k \phi(kx)$.
We can rewrite $g_{n,k}$ as
$$g_{n,k}(x,t)= \int_{0}^{t} \phi_{k}(B_{s}^{H,n}-x)ds.$$ 
We will now denote $I_{n,k}$ the discretisation of this integral calculated using the rectangle method using $m_{n,k}$ points uniformly spread on $[0,T]$
$$I_{n,k}(x,t)=\frac{1}{m_{n,k}} \sum_{i=0}^{[m_{n,k}t/T]} \phi_{k}\left( B_{\frac{i*T}{m_{n,k}}}^{H,n}-x \right),$$
where $[x]$ is the floor function.
We have $V_{n}(\cdot)=G_{n} e^{ 2 \abs{X_{n}}/ \alpha} l_{n}(X_{n},\cdot)$ that will be approximated by $W_{n,k}(\cdot)=G_{n} e^{ 2 \abs{X_{n}}/ \alpha} I_{n,k}(X_{n},\cdot).$ 

In the following, $C$ will denote a generic constant.

\begin{proposition}
Let $K$ denote a compact subset of $\R_{+}$. For $q>0$, for $\beta < 1/\alpha -1/2$, for all $\delta < \frac{1}{2H}-\frac{1}{2}$, if we take $m_{n,k}=[k^{\frac{\delta+2}{\delta'}} n^{-\frac{\beta}{\delta'}}]$ with $\delta'<H$,
there exists $C$ such that:
$$\E \left[\nr{V_{n}-W_{n,k}}_{K}^{q} \right] \leq \frac{Cn^{q\beta}}{k^{q\delta}}.$$
\end{proposition}

\begin{proof}
In this proof $C(\omega)$ will denote a generic random variable with finite moments of all order.

According to lemma \ref{lemmeholder},
\begin{equation}\label{diffvw}
\E \left[\nr{V_{n}-W_{n,k}}_{K}^{q} \right] \leq C \sup_{x\in \R}\left( \E\left[\nr{l_{n}(x,\cdot)-I_{n,k}(x,\cdot)}_{K}^{qp'}\right]\right)^{1/p'},
\end{equation}
Let us write
$$l_{n}(x,t)-I_{n,k}(x,t)= (l_{n}(x,t)-g_{n,k}(x,t)) + (g_{n,k}(x,t)-I_{n,k}(x,t)).$$
First, consider of $l_{n}(x,t)-g_{n,k}(x,t)=\int_{\R} \phi_{k}(y-x)(l_{n}(x,t)-l_{n}(y,t)) dy$. As the fractional
Brownian motion is locally non determinist, we can apply theorem 4 in \cite{pitt}, in order to have for all $\delta<\frac{1}{2H}-\frac{1}{2}$,
there exists $C(\omega)>0$ which has finite moments of all orders such that
\begin{equation}\label{difflg}
\abs{l_{n}(x,t)-g_{n,k}(x,t)} \leq \frac{C(\omega)}{k^{\delta}}.
\end{equation}
Now, let us consider 

\begin{align}\label{diffgi}
g_{n,k} & (x,t)-  I_{n,k}(x,t) =\nonumber\\
&\sum_{i=1}^{[m_{n,k}t/T]} \int_{(i-1)T/m_{n,k}}^{iT/m_{n,k}} \left( \phi_{k}(B_{s}^{H,n}-x) - \phi_{k} \left(B_{\frac{(i-1)T}{m_{n,k}}}^{H,n}-x \right) \right)ds
\nonumber\\
&+ \int_{ \left[ \frac{m_{n,k}t}{T}\right] \frac{T}{m_{n,k}}}^{t} \phi_{k}(B_{s}^{H,n}-x) ds - \phi_{k}\left(
B_{\left[\frac{m_{n,k}t}{T}\right]}^{H,n}-x\right).
\end{align}

We can remark that $\phi_{k}$ is $k^{2}$-Lipschitz, and that 
for $\delta'<H$ 
$$\sup_{s,t \in K} \frac{\abs{B_{t}^{H,n}-B_{s}^{H,n}}}{\abs{t-s}^{\delta'}}$$ 
has finite moments of all order. 
Consequently, using that majoration in (\ref{diffgi}) claims
\begin{equation}\label{diffgi2}
\abs{g_{n,k}(x,t)-I_{n,k}(x,t)} \leq \frac{k^{2}C(\omega)}{m_{n,k}^{\delta'}}.
\end{equation}

%
%
Combining (\ref{difflg}) and (\ref{diffgi2}) and by taking $m_{n,k}=[k^{\frac{\delta+2}{\delta'}} n^{-\frac{\beta}{\delta'}}]$
with $\beta < 1/\alpha -1/2$, there exists $C(\omega)$ with finite moments of all order such that
\begin{equation}\label{majorpsi}
\abs{l_{n}(x,t)-I_{n,k}(x,t)} \leq \frac{C(\omega)n^{\beta}}{k^{\delta}}.
\end{equation}
Using the majoration (\ref{majorpsi}) in (\ref{diffvw}) concludes.
\end{proof}

Using proposition \ref{evalsommelplq}, we get
\begin{corollary}\label{erreursomme}
For $p>0$ and $(N+1) \alpha >q > \max(p,2)$, for $\beta < 1/\alpha -1/2$, for all $\delta < \frac{1}{2H}-\frac{1}{2}$, if we take $m_{n,k}=[k^{\frac{\delta+2}{\delta'}} n^{-\frac{\beta}{\delta'}}]$ with $\delta'<H$,
there exists $C$ such that
$$\E \left[ \nr{\sum_{n=N+1}^{P}h(\Gamma_{n},V_{n}) - h(\Gamma_{n}, W_{n,k})}_{K,p}^{q} \right] \leq   \frac{C}{k^{q\delta}} \frac{1}{N^{q(2/\alpha-\beta-1)/2}}.$$
\end{corollary}

For the last part of this section, let us remind some notations
$$Y_{N}(t)= \sum_{n=1}^{N} \Gamma_{n}^{-1/ \alpha} G_{n} e^{ 2 \abs{X_{n}}/ \alpha} l_{n}(X_{n},t),$$
and $Y(t)=\lim_{N \rightarrow \infty} Y_{N}(t).$
Let us denote 
\begin{equation}\label{defz}
Z_{N,k} = \sum_{n=1}^{N} \Gamma_{n}^{-1/ \alpha} G_{n} e^{ 2 \abs{X_{n}}/ \alpha} I_{n,k}(X_{n},t).
\end{equation}

We will consider $\Pb(\nr{Y-Z_{P,k}}_{K,p}>3 \eps)$ for $(P+1)\alpha  > (N+1) \alpha >q > \max(p,2) \geq p \geq 1$.
\begin{align}
\Pb\left(\nr{Y-Z_{P,k}}_{K,p}> 3\eps\right) < &\ \Pb\left(\nr{Y-Y_{P}}_{K,p}>\eps\right) \nonumber\\
&+ \Pb\left(\nr{(Y_{P}-Y_{N})- (Z_{P,k}-Z_{N,k})}_{K,p}>\eps\right) \nonumber\\ 
&+ \Pb\left(\nr{Y_{N}-Z_{N,k}}_{K,p}>\eps\right). \nonumber
\end{align}

%

Corollaries \ref{erreurreste} and \ref{erreursomme} combined with Markov's inequality allow us to evaluate without any
difficulties both $$\Pb\left(\nr{Y-Y_{P}}_{K,p}>\eps \right) \leq \frac{C}{\eps^{q} P^{q(1/\alpha - 1/2)}}$$
 and $$\Pb\left(\nr{(Y_{P}-Y_{N})- (Z_{P,k}-Z_{N,k})}_{K,p}>\eps\right) < \frac{C}{\eps^{q}k^{q\delta}N^{q(1/\alpha-\beta-1/2)}}.$$

Now, we have to study the remaining terms $h(\Gamma_{n},V_{n}) - h(\Gamma_{n},W_{n,k})$ for $n \leq N$.

Denote $$\xi_{n,k}'(t)= \Gamma_{n}^{-1 / \alpha} G_{n} e^{2 \abs{X_{n}}/  \alpha} (l_{n}(x,t)-g_{n,k}(x,t)),$$
and 
$$\xi_{n,k}''(t)= \Gamma_{n}^{-1 / \alpha} G_{n} e^{2 \abs{X_{n}}/  \alpha} (g_{n,k}(x,t)-I_{n,k}(x,t)).$$
Since $\sum_{i=1}^{N} \xi_{n,k}'+\xi_{n,k}''= Y_{N}(t)-Z_{N,k}(t)$, we have
\begin{multline}
\Pb\left(\nr{Y_{N}-Z_{N,k}}_{K,p}>\eps \right) \leq \sum_{n=1}^{N} \Pb\left(\nr{\xi_{n,k}'}_{K,p}> \frac{\eps}{2N}\right)\\
+\sum_{n=1}^{N} \Pb\left(\nr{\xi_{n,k}''}_{K,p}>\frac{\eps}{2N}\right).
\end{multline}


A direct computation is enough to evaluate $$\Pb\left(\nr{\xi_{n,k}'}_{K,p}> \frac{\eps}{2N}\right) \leq \frac{(2N)^{\alpha}C}{\eps^{\alpha}k^{\alpha
\delta}}.$$
 
For $n \leq N$ taking $m_{n,k}=[\Gamma_{n}^{-1/(\delta' \alpha)}k^{\frac{2+\delta}{\delta'}}]$ we have $$\E \left[ \frac{\Gamma_{n}^{-q/\alpha}}{m_{n,k}^{\delta' q}}\right]
\leq\frac{C}{k^{q(\delta+2)}},$$ so we can apply Markov's inequality on the last term. (\ref{diffgi2}) leads to the
evalution of $$\Pb\left(\nr{\xi_{n,k}''}_{K,p}>\frac{\eps}{2N}\right) \leq \frac{(2N)^{q}C}{\eps^{q}k^{q\delta}}.$$ 

Once this remaks have been made, we only have to tune our parameters to control the error before stating the following
theorem:
\begin{proc}\label{assfinal}
We want to approximate process $Y$ of parameters $(H,\alpha)$ by a family $(Z_{P_{\eps},k_{\eps}})_{\eps >0}$.
We see process $Y$ as a series and $Z_{P_{\eps},k_{\eps}}$ as a truncated series.

We can adjust two parameters: $P_{\eps}$ is the size of the truncature and $k_{\eps}$ controls the approximation of the fractional Brownian motion local time approximation.
We will make two kinds of errors, one coming from the truncature itself, and one from the approximation of the local time.

The error from the truncature can be controlled if we take $P_{\eps} \sim C \eps^{-\eta\frac{2\alpha}{2-\alpha}}$.

The approximation error has two sources: our theorical approximation of the local time, and the discretisation used to compute the approximation.

We can deal with the first one by taking  $k_{\eps} \sim C \eps^{-\eta/\delta}$ where $\delta<1/(2H)-1/2$ comes from the Hölder continuity of the fractional Brownian motion local time.

Let us denote by $m_{n,k}$ the number of points used in the discretisation. 
In the series defining $Y$, the first terms are the most important so we will distinguish two cases, the first $N$ terms, where $p>1$ and $(P+1)\alpha  > (N+1) \alpha >q > \max(p,2)$, where we will take  $m_{n,k}$ random, linked with the size of $\Gamma_{n}$, $m_{n,k}=[\Gamma_{n}^{-1/(\delta' \alpha)}k_{\eps}^{\frac{2+\delta}{\delta'}}]$, where $\delta'<1/H$ comes from the Hölder continuity of the fractional Brownian motion.
For the remaining terms, a good precision is not as important, so we will need fewer points in our discretisation. We can take $m_{n,k}=[k^{\frac{\delta+2}{\delta'}} n^{-\frac{\beta}{\delta'}}]$ with $\beta<1/\alpha -1/2$.
\end{proc}

\begin{theorem}\label{thmfinal}
With tuning procedure \ref{assfinal}, we are able to get a convergence rate for family of process $Z_{P_{\eps},k_{\eps}}$ defined in (\ref{defz}). For $\eta > 1$ and $p\geq1$ there exists $C$ such that
$$P(\nr{Y-Z_{P_{\eps},k_{\eps}}}_{K,p} > \eps) \leq C \eps^{\alpha (\eta - 1)}.$$
\end{theorem}

Using $P_{\eps}$, $k_{\eps}$ and $m_{n,k}$ as defined previously, we can have one last theorem:
\begin{theorem}
$Z_{P_{\eps},k_{\eps}}$ converges almost surely towards $Y$.
\end{theorem}

\begin{proof}

For the assumptions used in theorem \ref{thmfinal} with $q \geq 2$ and  $N>q / \alpha-1$ fixed.

Using Minkowski's inequality: 
\begin{align}\label{ineqtriang}
\E \Bigl[ \nr{(Y(t)-Y_{N}(t))-\right. & \left. (Z_{P_{\eps},k_{\eps}}-Z_{N,k_{\eps}})}_{K,p}^{q}   \Bigr]^{1/q}\leq\nonumber\\
&\E \left[ \nr{(Y-Y_{P_{\eps}})}_{K,p}^{q}   \right]^{1/q} \nonumber\\
+ &\E \left[ \nr{(Y_{P_{\eps}}-Y_{N})- (Z_{P_{\eps},k_{\eps}}-Z_{N,k_{\eps}})}_{K,p}^{q}   \right]^{1/q}.     \nonumber
\end{align}

Thanks to corollaries \ref{erreurreste} and \ref{erreursomme}, using the expression of $P_{\eps}$, $k_{\eps}$ and $m_{n,k}$
given in \ref{assfinal}, we get each term bounded by $C \eps^{\eta}$, with $\eta>1$. Thus, we can say that
\begin{equation}\label{erreurreste2}
\E \left[ \nr{(Y-Y_{N})- (Z_{P_{\eps},k_{\eps}}-Z_{N,k_{\eps}})}_{K,p}^{q}   \right]^{1/q} \leq C \eps^{q \eta}
\end{equation}

Denote $q$ such as $q \eta>1$. Inequality (\ref{erreurreste2}) and Borel-Cantelli lemma imply $Z_{P_{\eps},k_{\eps}}-Z_{N,k_{\eps}}$ converges towards $Y-Y_{N}$.
We only have to prove $Z_{N,k_{\eps}}$ converges almost surely towards $Y_{N}$. But according to (\ref{majorpsi}), using $k_{\eps}$ and $m_{n,k}$ chosen, $\nr{l_{n}(x,t)-I_{n,k}(x,t)}_{K}$ converges almost surely towards $0$ so $Y_{N}-Z_{N,k_{\eps}}$ too since only a finite number of terms is considered.
\end{proof}

\section{Simulation}

\begin{figure}
\begin{center}
\includegraphics[height=70mm]{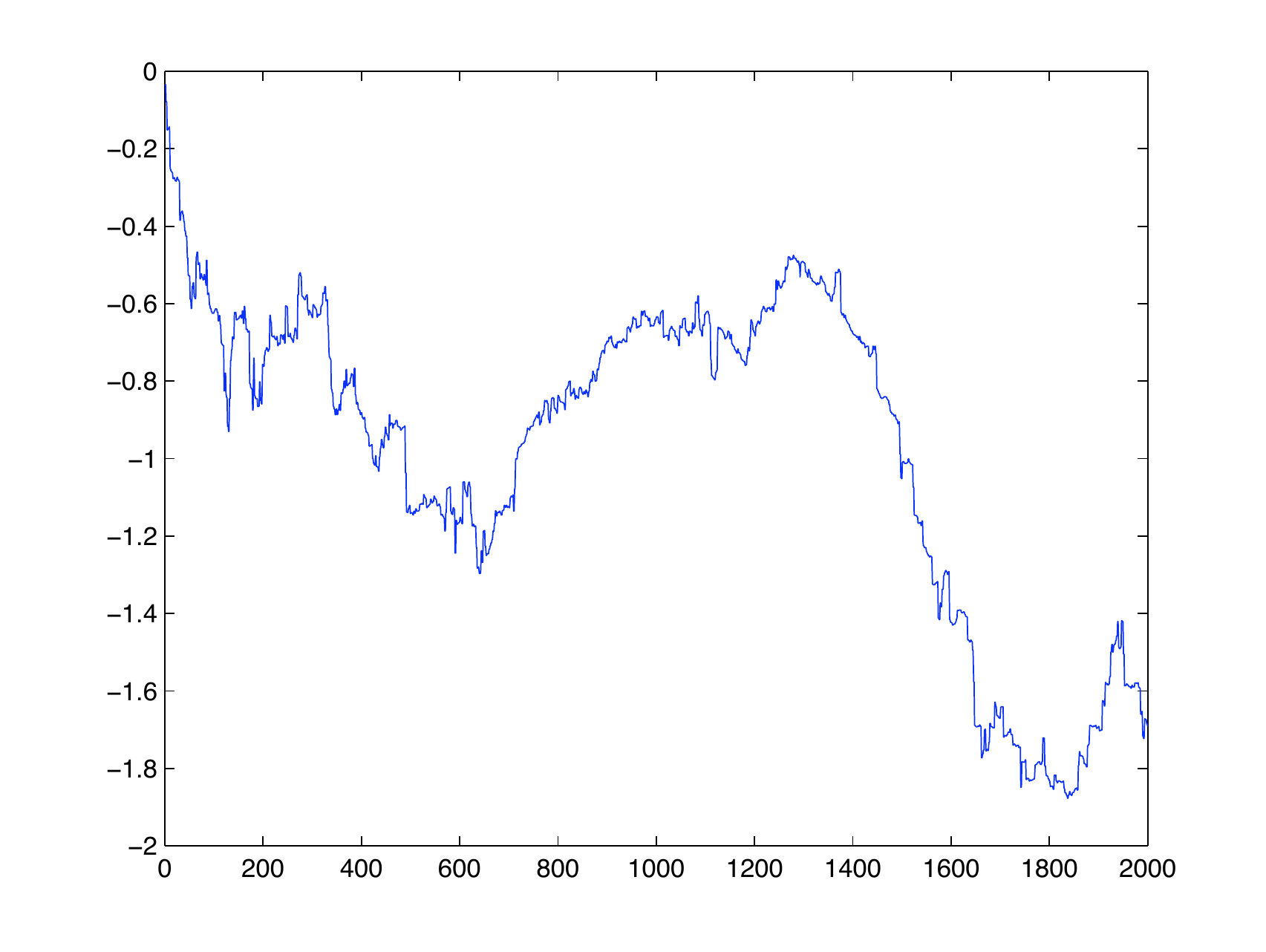}
\caption{Allure of $Y$ for $H=0.3$ and $\alpha=1.2$}
\label{samplepath}
\end{center}
\end{figure}

In \cite{cohensamo}, the authors explained how to simulate the Brownian motion local time stable motion with a random walk with random rewards approach. They are many advantages of our approach against the random walk with random reawards approach. The most obvious is that it is valid for all Hurst parameter and not only $H=1/2$. Moreover we already highlighted that we have a convergence rate which was not the case previously.

Unfortunately, as one can see in our tuning procedure, the parameters we use depend on $\alpha$ and $H$. If $\alpha$ is close to $2$, the number of terms in the sum is to high to be accepted and if $H$ is close to $1$, the precision needed in the simulation of the fractional Brownian motion is also to high to do. Even if this method is not perfect for all values of $H$ and $\alpha$, it is still a major improvement since we can have sample paths for many different values of $H$ with a convergence rate. See for example one simulation in figure \ref{samplepath}.

\begin{figure}
\begin{center}
\includegraphics[height=60mm]{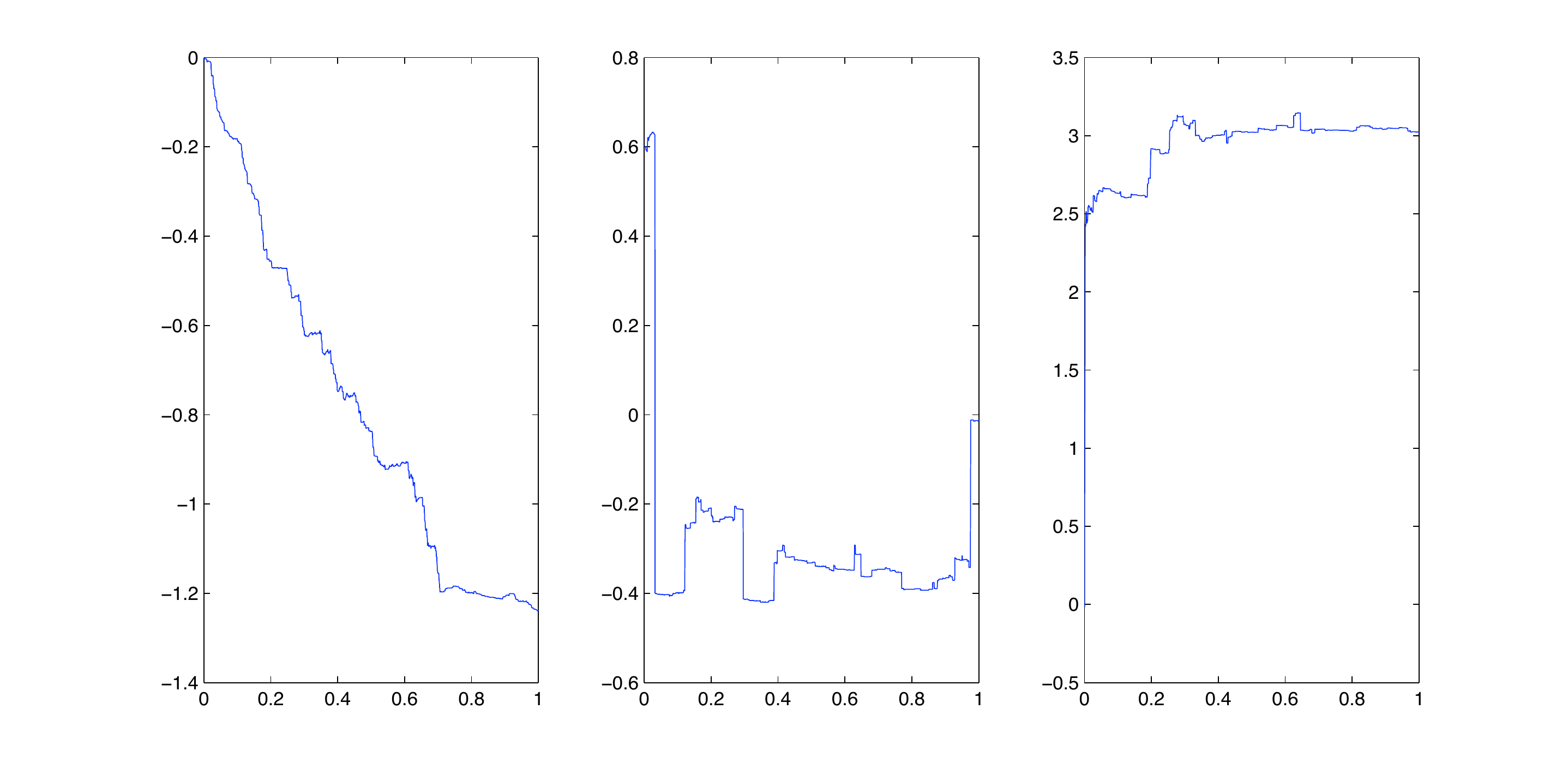}
\caption{Allure of $Y$ for $H=0.2$,  $H=0.4$ and $H=0.6$ with $\alpha=0.7$}
\label{samplepaths}
\end{center}
\end{figure}

According to theorem 5.1 in \cite{cohensamo}, the Hölder exponent $d$ of our process is such that $d<1-H$.
It means that the closer to $1$ is $H$, the less regular is our process. See figure \ref{samplepaths} for sample paths with different values of $H$ and $\alpha$ constant.


\begin{figure}
\begin{center}
\includegraphics[height=70mm]{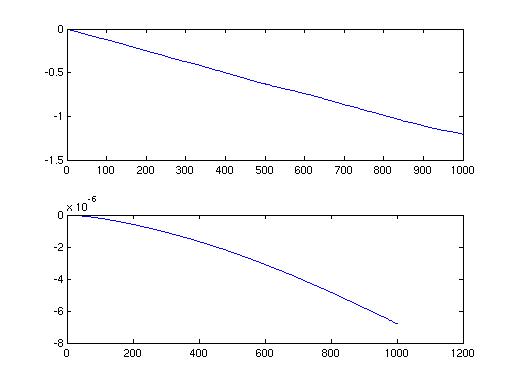}
\end{center}
\caption{
Allures of  $\log (\E \left[e^{i u Y(t)} \right])$ with respect to $t$. Our method on the first line, and the random walk method on the second.
}
\label{droites}
\end{figure}
%
%
%

Let us try to make a more rigorous comparison of our two ways to simulate our process. We have to take $H=0.5$, so that both methods are able to do the simulation.
In order to check the validity of both approaches, we are going to compute $\E \left[e^{i u Y(t)} \right]$. A straightforward computation from the result of proposition \ref{carac}, gives in the case $\alpha=1$,
$$\E \left[e^{i u Y(t)} \right] = \exp \left(C \abs{u} t \right).$$
We are going to check for both of our processes if $\log (\E \left[e^{i u Y(t)} \right])$ is a straight line using a classical Monte Carlo method. 
The result is shown in figure \ref{droites}, both methods using $10 000$ steps, the first one being the method developed in this paper, the second the random walk with random rewards method.

According to a linear regression, the $R^{2}$ coefficient is $0.9997$ for the first one, and $R^{2}=0.9763$.
Thus we can conclude that even if the random walk method seems a bit quicker, our method has not only the advantage of being able to deal with $H \neq 0.5$ but also is closer to what we should theoretically get.

\bibliography{ltfsm}

\end{document}